\documentclass[twoside]{aiml22}

\usepackage{aiml22macro}

\usepackage{graphicx}
\usepackage{amsmath}
\usepackage{amssymb}
\usepackage{comment}
\usepackage{tikz-cd}

\newcommand{\op}{\overline{p}}
\newcommand{\oq}{\overline{q}}
\newcommand{\orr}{\overline{r}}
\newcommand{\opsi}{\overline{\psi}}

\def\lastname{Almeida and Ghilardi}

\begin{document}

\begin{frontmatter}
  \title{Unification with Simple Variable Restrictions and Admissibility of $\Pi_{2}$-Rules}
  \author{Rodrigo Nicolau Almeida}\footnote{r.dacruzsilvapinadealmeida@uva.nl.}
  \address{Institute for Logic Language and Computation (ILLC) - University of Amsterdam \\ Science Park 107 \\ 1098 XG Amsterdam}
  \author{Silvio Ghilardi}\footnote{silvio.ghilardi@unimi.it}
  \address{{Dipartimento di Matematica} \\ {via C. saldini 50} \\ {20133 Milano}}

  \begin{abstract}
  We develop a method to recognize admissibility of $\Pi_{2}$-rules, relating this problem to a specific instance of the unification problem with linear constants restriction \cite{BAADER1996211}, called here ``unification with simple variable restriction".
  It is shown that for logical systems enjoying an appropriate algebraic semantics and a 
  finite approximation of left uniform interpolation,
this unification with simple variable restriction can be reduced to standard unification.
As a corollary, we obtain the decidability of admissibility of $\Pi_{2}$-rules for many logical systems.
  \end{abstract}

 \end{frontmatter}

\section{Introduction}\label{sec:intro}

Non-standard rules have often been used in the context of logical systems to axiomatise specific classes of models. Their use traces its origin to the work of Gabbay \cite{Gabbay1981} as well as Takeuti and Titani \cite{Takeuti1984-TAKIFL}, and has been the subject of some attention, especially with a focus on axiomatisation and admissibility of such rules \cite{venemaantiaxioms,BEZHANISHVILI-stricimplicationcalculus,Bezhanishvili2022-if,rodrigopaperpi2rules}. Nevertheless, in most of these contexts specific assumptions have been made on what counts as such a non-standard rule, which make it difficult to provide a unified account of what these rules should be, and which make the current results available in the literature difficult to transfer: for instance, whilst in \cite{Bezhanishvili2022-if}, some connections were made between the solution of the admissibility problem for some modal logic systems, the existence of uniform interpolants, and some problems of unification, it is not clear how to generalize this to settings such as the Takeuti and Titani rule.

In this paper we start from the simple observation that, when left uniform interpolants are available, admissibility of $\Pi_2$ rules can be reduced to admissibility of standard rules just by eliminating bound context  variables via such interpolants. It is less obvious that one can get the same result when assuming only that left uniform interpolants are 'finitely approximable': in fact, in this case  one needs to show that such finite approximations are stable under substitutions. We obtain the result by employing techniques from two different sources: on one side, we reduce our task to subobject manipulations in the opposite category of finitely presented algebras (in the style of~\cite{Ghilardizawadowskybook}) and on the other side we connect admissibility problems for $\Pi_2$ rules to a dedicated
 version of $E$-unification theory obtained by 
 specializing the  
 ``Unification with Linear Constant Restrictions" employed in
  \cite{BAADER1996211} to handle combined $E$-unificaton problems and general $E$-unification problems.


The structure of the paper is as follows: in Section \ref{Admissibility of Pi2-rules} we define
formally the problem of admissibility of $\Pi_{2}$-rules, and state our main result. 
In Section \ref{Unification with Simple Constant Restrictions} we introduce the problem of unification with simple variable restriction, and provide an equivalent algebraic  presentation of it. In Section \ref{Algebraic Correspondents to Interpolation} we recall the correspondence between some logical properties we  need and their reformulations in the opposite of the category of finitely presented algebras.
 In Section \ref{Unification in Exceptional Logical Systems} we prove our main theorem, showing that under suitable assumptions, the unification type for simple variable restriction is finitary and reduces to the standard unification type.  In Section \ref{Applications and Negative Results} we provide some applications;  in Section \ref{Conclusions and Further Work} we conclude and highlight some limitations  of our method.
 In Section~\ref{NIS}, we analyze the prominent example of nuclear implicative semilattices.

\section{Admissibility of $\Pi_{2}$-rules in Logical Systems}\label{Admissibility of Pi2-rules}

Throughout we will assume that we are working in a functional signature $\fancyL$
comprising at least a constant symbol; 
the set of \emph{terms} (aka propositional formulas, or just \emph{formulas})  is denoted by $Fm_{\fancyL}$. A \emph{logic} $\vdash$ in this language is a relation
$\vdash\;\subseteq \wp(Fm_{\fancyL}) \times Fm_{\fancyL}$ satisfying the usual identity, monotonicity,  transitivity, structurality (i.e. invariance under substitutions) and finitarity conditions
(see \cite[Definitions 1.5-1.6]{Font2016-dk}). We use the letters $\phi, \psi,\dots$ for $\mathcal{L}$-formulas and letters $p, q,\dots$  or $x,y, \dots$ for variables; we compactly represent a tuple of distinct variables
as $\op$. The notation $\phi(\op)$ means that the formula $\phi$ has free variables included in the tuple $\op$. Since our \emph{tuples of variables} are assumed to be formed by \emph{distinct} elements,  we emphasize that when we write e.g. $\phi(\op, \oq)$, we mean that the tuples $\op, \oq$ are made of distinct variables and are also disjoint from each other. Notations like $\phi(\opsi/\op)$ (or just $\phi(\opsi)$)
denote the result of substituting $\overline{p}$ by $\overline{\psi}$ inside of $\phi$.
If $\Gamma, \Delta$ are sets of formulas, $\Gamma\vdash \Delta$ means $(\Gamma,\phi)\in\; \vdash$ for all $\phi\in \Delta$.

When a logic is algebraizable, most of the definitions we shall introduce  can be transferred back and forth from the corresponding class of algebras. Recall that $\vdash$ is algebraizable  iff there is a quasivariety of $\fancyL$-algebras  $\bf K$ and there are essentially inverse structural transformers between $\fancyL$-formulas and $\fancyL$-equations, mapping elements of $\vdash$ to $\bf K$-valid quasi-equations and vice versa
    (see~\cite[Definitions 3.11]{Font2016-dk}).
    We say that our  logic $\vdash$ (which we assumed to be finitary)  is \emph{strongly algebraizable} iff $\mathbf{K}$ is actually a variety.
    
\begin{assum}
 For the whole paper, we fix a language $\fancyL$ and a strongly algebraizable logic $\vdash$ in it; we call $\bf K$ the equivalent algebraic semantics of $\fancyL$-algebras and $E$ an equational theory axiomatizing $\bf K$. 
\end{assum}

We shall also fix a Hilbert-style derivation system $\vdash_S$ associated to the logic $\vdash$; vacuously one always exists (simply considering the set of rules $\Gamma\vdash \phi$ whenever $\Gamma$ is finite and $\Gamma\vdash\phi$, see \cite{Font2016-dk} for further details).
%
%
We will begin by outlining in general what a $\Pi_{2}$-rule is in this context; the definitions here are analogous to the ones presented in \cite{rodrigopaperpi2rules}.

\begin{definition}
Let $\Gamma=\{\phi_{i}(\overline{p},\overline{q}): i\leq n \}$ and $\psi(\overline{q})$ be formulas in the language $\fancyL$.
The $\Pi_{2}$-\textit{rule} associated with this sequence of formulas is denoted $\forall\overline{p}\Gamma/^{2}\psi$ (sometimes without the universal quantifier, when the variables are clear from context)\footnote{This notation serves to emphasise
both the fact that
these are distinct from usual rules, and the second-order nature of these rules, but we point out that it is purely formal.} and usually
displayed as:
\begin{prooftree}
\AxiomC{$
\forall \overline{p}~(
\phi_{0}(\overline{p},\overline{q}),...,\phi_{n}(\overline{p},\overline{q}))
$}
\UnaryInfC{$\psi(\overline{q}).$}
\end{prooftree}
Given such a collection of formulas $\Gamma=\{\phi_{i}\}_i$, we refer to
$\overline{p}$ as the \textit{bound context of $\Gamma$}, or generally, the \textit{bound context associated to $\Gamma$}, and sometimes denote it as $F_{c}(\phi_{i}$) or $F_{c}(\Gamma)$; we refer to propositional variables not ocurring in $F_{c}$ as the \textit{free context}. Whenever the bound context is empty, the rule is referred to as a \emph{standard} rule.
\end{definition}

\begin{example}
    Let $\fancyL$ be the language of modal logic. Gabbay's \textit{irreflexivity rule} is the rule
    \begin{prooftree}
        \AxiomC{$\forall p\,((\Box p\rightarrow p)\vee \phi)$}
        \UnaryInfC{$\phi$}
    \end{prooftree}
    where $\phi$ is any formula, such that $p$ does not occur in $\phi$. This was used in \cite{Gabbay1981} to obtain completeness with respect to a class of irreflexive frames.
\end{example}

\begin{example}
    Let $\fancyL$ be the language of modal algebras with a binary modality $\rightsquigarrow$, called the signature of \textit{contact algebras}. Consider the following rule:
    \begin{prooftree}
        \AxiomC{$\forall p\,((p\rightsquigarrow p)\wedge (\phi\rightsquigarrow p)\wedge (p\rightsquigarrow \psi)) \rightarrow \chi$}
        \UnaryInfC{$(\phi\rightsquigarrow \psi)\rightarrow \chi.$}
    \end{prooftree}
    This rule was discussed in \cite{BEZHANISHVILI-stricimplicationcalculus} and \cite{Bezhanishvili2022-if}, and used to axiomatise the strict implication calculus.
    Notice however that our notion of $\Pi_2$-rule is more general than the $\Pi_2$-rules introduced in \cite{BEZHANISHVILI-stricimplicationcalculus} and \cite{Bezhanishvili2022-if}
    (for instance, the $\Pi_2$-rules introduced there do not have all standard rules as special cases).
\end{example}

\begin{example}\label{Takeuti-Titani Rule}
    Let $\fancyL$ be the language of Heyting algebras. Consider the following rule, often called the \textit{Takeuti-Titani rule} or the \textit{density rule}:
    \begin{prooftree}
        \AxiomC{$\forall r(g\rightarrow ((p\rightarrow r)\vee (r\rightarrow q)\vee c))$}
        \UnaryInfC{$g\rightarrow (p\rightarrow q) \vee c$}
    \end{prooftree}
    Such a rule has been fruitfully used to axiomatise classes of G\"{o}del algebras and other MV-algebras (see e.g. \cite{metcalfemontagnasubstructuralfuzzy,Baaz2017}). We will return to it as an example later.
\end{example}

 We now explain how $\Pi_{2}$-rules can be used within the derivation system $\vdash_S$:

\begin{definition}\label{Derivation using pi2-rules}
Let $\Sigma$ be a set of $\Pi_{2}$-rules. Given a formula $\phi$ we say that $\phi$ is \textit{derivable} using the $\Pi_{2}$-rules in $\Sigma$, and write $\vdash_{S\oplus \Sigma}\phi$, provided there is a sequence $\psi_{0},...,\psi_{n}$ of formulas such that:
\begin{itemize}
\item $\psi_{n}=\phi$;
\item For each $\psi_{i}$ we have that either:
\begin{enumerate}
\item $\psi_{i}$ is an instance of an axiom of $\vdash_{S}$ or,
\item $\psi_{i}$ is obtained using a rule from $\vdash_{S}$, from some previous $\psi_{j_{0}},...,\psi_{j_{k}}$ or,
\item $\psi_{i}= \chi(\overline{\xi}/\oq)$ and $\psi_{j_{k}}= \mu_{k}(\overline{r}/\op,\overline{\xi}/\oq)$ for $0\leq j_{k}<i\leq n$, where
\begin{itemize}
    \item  $\overline{r}$ is a renaming of $\op$, away from $\overline{\xi}$, i.e., a set of fresh variables not ocurring in $\overline{\xi}$;
    \item $\Delta=\{\mu_{k}(\overline{p},\overline{q}) : k\in \{0,...,m\}\}$;
    \item $\forall\op\Delta/^{2}\chi\in \Sigma$;
    \item $\chi=\chi(\overline{q})$.
\end{itemize}
\end{enumerate}
\end{itemize}
\end{definition}

An \emph{extended calculus} is a calculus of the kind $\vdash_{S\oplus \Sigma}$.
We can now write what it means for a rule to be \textit{admissible}:

\begin{definition}
    Let $\Gamma/^{2}\phi$ be a $\Pi_{2}$-rule, and $\vdash_{S\oplus\Sigma}$ some extended calculus. We say that the rule $\Gamma/^{2}\phi$ is \textit{admissible} in $\vdash_{S\oplus \Sigma}$ if for all $\psi$:
    \begin{equation*}
        \vdash_{S\oplus\Sigma\oplus \Gamma/^{2}\phi}\psi \implies \vdash_{S\oplus \Sigma}\psi.
    \end{equation*}
\end{definition}

In light of the  definition of derivation using non-standard rules,
we want to ``internalize'' the notion of admissibility for ordinary (i.e. non extended) calculi.
A standard rule is admissible over $\vdash$ iff in every substitution making the premises into theorems of $\vdash$ also makes the conclusion into a theorem of $\vdash$.\footnote{It should be noticed however that this characterization does not always hold, for instance it fails for multiple-conclusion rules, see~\cite{MetcalfeWolliC}
for a thorough discussion.} In order to obtain a similar characterization here, we need the notion of $C$-invariant substitution. Given a finite set $C=\{p_{1},...,p_{n}\}$ of propositional variables,
a \emph{$C$-invariant substitution} is a substitution $\sigma$ mapping the $p_i$ into themselves and the other propositional variables $q$ into formulas $\sigma(q)$ not containing the variables in $C$.

\begin{lemma}\label{Internal Characterization of Admissibility}
    Let $\forall\overline{p}\Gamma/^{2}\phi$ be a $\Pi_{2}$-rule and let $C:=\{p_{1},...,p_{n}\}$. Then $\forall\overline{p}\Gamma/^{2}\phi$ is admissible over $\vdash_{S}$ if and only if whenever $\sigma$ is a $C$-invariant substitution and we have $\vdash_{S}\sigma(\Gamma)$, then we have also $\vdash_{S}\sigma(\phi)$.
\end{lemma}
\begin{proof}
    Follows by induction on the structure of the derivation
    and by the fact that $\vdash_{S}$ is invariant under variable renamings.
\end{proof}

Given a logic $\vdash$,
by the \textit{$\vdash$-admissibility problem for $\Pi_{2}$-rules}  we mean the problem of  determining, given a triple $(\Gamma,\phi,C)$ (where $C=\overline{p}$ is a set of propositional letters ocurring in $\Gamma$
but not in $\phi$),
whether the $\Pi_{2}$-rule $\forall \overline{p}\Gamma/^{2}\phi$ is admissible over $\vdash$. It is well-known that the $\vdash$-admissibility problem for
standard rules need not be decidable (see for example \cite[Theorem 4]{Baader2010}), and that it is decidable for a substantial number of logical systems encountered in practice, such as $\mathsf{S4},\mathsf{S5},\mathsf{IPC}$ and lax logic~\cite{GhLe} amongst others. In this paper we will obtain a general result concerning the decidability of admissibility for $\Pi_{2}$-rules, under specific assumptions which we now proceed to review.

First, we need to recall some notions concerning different interpolation properties.
The following definition is modelled after \cite[Section 2]{Czelakowski1985}:

\begin{definition}\label{def:mint}
We say that a logic $\vdash$ has the \textit{Maehara interpolation property} if for any finite sets $\overline{p},\overline{q},\orr$ of propositional variables and for any finite  sets of formulas $\Sigma(\op,\oq),\Delta(\oq,\orr), \Sigma'(\op,\oq)$ in the language $\fancyL$  the following holds: if $\Sigma(\op,\oq)\cup \Delta(\oq,\orr)\vdash \Sigma'(\op,\oq)$, then there exists a set of formulas $\Pi(\oq)$, such that  $\Delta(\oq,\orr)\vdash \Pi(\oq)$ and $\Sigma(\op,\oq)\cup \Pi(\oq)\vdash \Sigma'(\op,\oq)$.
\end{definition}

\begin{definition}\label{def:rint}
 We say that a logic $\vdash$ has \textit{right uniform deductive interpolation} if for any finite sets $\overline{p},\overline{q}$ of propositional variables, and finite set of formulas $\Sigma(\overline{p},\overline{q})$ there exists a finite set of formulas $\Pi(\overline{q})$ such that for any finite set of propositional variables $\orr$ and for any finite set of formulas $\Delta(\overline{q},\overline{r})$
    \begin{equation*}
        \Sigma(\op,\oq)\vdash \Delta(\oq,\orr) \iff \Pi(\oq)\vdash \Delta(\oq,\orr).
    \end{equation*}
\end{definition}

    There is a specular
    \textit{left uniform deductive interpolation property} saying that
    for any finite sets $\overline{q},\overline{r}$ of propositional variables, and finite set of formulas $\Delta(\overline{q},\overline{r})$ there exists a finite set of formulas $\Theta(\overline{q})$ such that for any finite set of propositional variables $\op$ and for any set of formulas $\Sigma(\overline{p},\overline{q})$
    \begin{equation*}
        \Sigma(\op,\oq)\vdash \Delta(\oq,\orr) \iff \Sigma(\op,\oq)\vdash \Theta(\oq).
    \end{equation*}
    %
    %
    The two properties are not equivalent \footnote{See
    Example~\ref{ex:isl} below
    for a counterexample. Equivalence can hold in very special contexts, typically when compact congruences are Boolean (this is the case of~\cite[Section 4]{Bezhanishvili2022-if} for the presence of a universal modality).}; the latter property is sometimes denoted by saying that the logic $\vdash$ has \textit{global post-interpolants}. In this paper, we consider  the following strictly weaker  version of left uniform deductive interpolation:

    \begin{definition}\label{def:lint}
    We say that $\vdash$ has \textit{left-finitary uniform deductive interpolation} if
    for any finite sets $\overline{q},\overline{r}$ of propositional variables, and finite set of formulas $\Delta(\overline{q},\overline{r})$
     there is a finite collection of finite sets of formulas $\Theta_{1}(\overline{q}), \dots,\Theta_{n}(\overline{q}) $ such that:
    \begin{align*}
        &{\text (i)}~\Theta_{i}(\oq)\vdash \Delta(\oq,\orr) \text{ for each $i\leq n$ and } \\
        &{\text (ii)} \text{ for any $\op$ and $\Sigma(\overline{p},\overline{q})$, } \\
        & ~~~~~~~~~~\Sigma(\op,\oq)\vdash \Delta(\oq,\orr) \implies \Sigma(\op,\oq)\vdash\Theta_{i}(\oq)
        \\
        & ~~~~\text{ for some $i\leq n$. }
    \end{align*}
    \end{definition}

When we move to concrete decidability problems and say that $\vdash$ has left-finitary uniform deductive interpolation, we also assume that the above finite sets $\Theta_{1}(\overline{q}), \dots,\Theta_{n}(\overline{q}) $ are \emph{effectively computable} from $\Delta(\overline{q},\overline{r})$.
The main result of the paper is as follows:
    
\begin{theorem}\label{Key Theorem}
    Suppose that $\vdash$ has
    \begin{enumerate}
        \item the Maehara Interpolation Property;
        \item Right-Uniform Deductive Interpolation;
        \item Left-Finitary Deductive Uniform Interpolation.
    \end{enumerate}
    Then if both $\vdash$ itself and the  $\vdash$-admissibility problem for standard rules  are decidable, so is the $\vdash$-admissibility problem for $\Pi_{2}$-rules.
\end{theorem}

After establishing the algebraic analogues of these syntactic properties (in Section 4), the proof of this Theorem will be given as Corollary \ref{Corollary connecting unification types and admissibility problems}.
As a consequence, we obtain that several well-studied logical systems have such a decidable problem: among them we have
%
$\mathsf{S5},\mathsf{GL},\mathsf{Grz}, \mathsf{IPC},  \mathsf{LC} $ and
lax logic~\cite{rosalie}, all of which satisfy the hypotheses of the Theorem. As for logics without left uniform deductive interpolation property, we mention
the $\{\wedge,\top, \rightarrow\}$ -fragment of $\mathsf{IPC} $ and the
$\{\ell,\wedge,\top, \rightarrow\}$-fragment of lax logic: these systems  have Maehara interpolation, are locally finite (by Diego theorem and extensions~\cite{BBCGGJ})
and have decidable admissibility problems for standard rules.
In Section~\ref{Applications and Negative Results}, we illustrate with examples how the algorithm from Theorem~\ref{Key Theorem} works.


\section{Unification with Simple Variable Restrictions}\label{Unification with Simple Constant Restrictions}

In this section we recall some essential concepts from unification theory which will be needed in our work; the reader can find more general information on unification and $E$-unification theory in \cite{siekmanbaaderingabbay} and the references contained therein.
Our aim is to introduce \textit{unification with simple variable restrictions}:
this is a special case of
 unification with linear constant restrictions as introduced in \cite{BAADER1996211}, where
 unification with linear constant restrictions is an essential ingredient for building combined unification algorithms.

We let $Fm_{\fancyL}(\op)$ be the set of ${\fancyL}$-formulas containing at most the variables $\op$; in algebraic terms, $Fm_{\fancyL}(\op)$ is the absolutely free  algebra over $\op$ (also called the \emph{term algebra} over $\op$).
A \textit{substitution}
is an ${\fancyL}$-morphism of  term algebras $\sigma:Fm_{\fancyL}(\op)\to Fm_{\fancyL}(\oq)$;
hence $\sigma$ can be represented as a finite set of variable-term pairs,
\begin{equation*}
    \sigma=\{p_{1}\leftarrow \phi_{1}(\oq),...,p_{n}\leftarrow \phi_{n}(\oq)\}.
\end{equation*}
We say that $\op$ is the domain and $\oq$ the co-domain of the substitution $\sigma$.
Recall that  a substitution $\sigma$ is said to be \emph{$C$-invariant}
(for a finite set $C$ of propositional variables included into the domain and in the codomain of $\sigma$) iff it maps the $p\in C$
  into themselves and the $q\not \in C$ into terms $\sigma(q)$ not containing the variables in $C$.

Given
two terms $\phi,\psi$, and an equational theory $E$, we write $\phi=_{E}\psi$ to mean that
$E \models \phi=\psi$.\footnote{Of course, since we assume algebraizability, all definitions in this section could be equivalent stated inside the logical context of $\vdash$ by applying the appropriate transformers.}
We say that two substitutions $\sigma$ and $\tau$ having the same domain and codomain are $E$-equivalent, briefly $\sigma=_{E}\tau$, if and only if $\sigma(p)=_{E}\tau(p)$ holds for every variable $p$ in their domain. We say that $\tau:Fm_{\fancyL}(\op)\to Fm_{\fancyL}(\oq)$ is \textit{less general than $\sigma:Fm_{\fancyL}(\op)\to Fm_{\fancyL}(\oq')$} (with respect to $E$) if there is a substitution $\theta:Fm_{\fancyL}(\oq)\to Fm_{\fancyL}(\oq')$ such that
\begin{equation*}
    \tau=_{E}\theta\circ \sigma.
\end{equation*}

\begin{definition} Given a finite set of propositional variables $C$,
an $E$-\textit{unification problem with simple variable restriction}
(briefly a \emph{C-unification problem}) is a finite set
of pairs of terms in the variables $\op$ (with $C\subseteq \op)$\footnote{
In this paper, we do not consider free constants -- i.e., fixed propositional variables which are interpreted freely in our algebras -- in unification problems; considering them, would lead to consider \emph{parameters} in inference rules.}
\begin{equation*}
    (\phi_{1}(\op),\psi_{1}(\op)),...,(\phi_{k}(\op),\psi_{k}(\op)); \leqno (P_C)
\end{equation*}
a solution to such a problem or a \textit{$C$-unifier} is a $C$-invariant substitution $\sigma$ of domain $Fm_{\fancyL}(\op)$ such that
\begin{equation*}
    \sigma(\phi_{1})=_{E}\sigma(\psi_{1}),...,\sigma(\phi_{k})=_{E}\sigma(\psi_{k}).
\end{equation*}
\end{definition}
When $C=\emptyset$, we speak of \emph{standard} unification problems, or just unification problems. We also write $U_{E}^{svr}(P_C)$ for the set of $C$-unifiers for the problem $P_{C}$.


Once  $C$ and  $(P_C)$ (a unification problem as above) are fixed, given two $C$-unifiers $\tau$ and $\sigma$, we say that $\tau$ is \textit{less general} than $\sigma$, and write $\tau\leq_{C}\sigma$ if it is less general
than $\sigma$
as a substitution. Hence, given a unification problem $(P_C)$ with simple variable restrictions, this definition of $\leq_C$ induces a preorder on $U_{E}^{svr}(P_C)$. A \emph{$C$-unification basis} from this set is a subset $B\subseteq U_{E}^{svr}(P_C)$ such that for every $\sigma'\in U_{E}^{svr}(P_C)$ there is $\sigma\in B$ such that $\sigma'\leq_{C}\sigma$ holds; a \textit{most general $C$-unifier} ($C$-mgu) of $(P_C)$
is a $\sigma\in U_{E}^{svr}(P_C)$ such that $\{\sigma\}$ is a $C$-unification basis.


\begin{definition}\label{def:utype}
    We say that $E$ has \textit{ finitary simple-variable-restriction (svr) unification type}  iff every $C$-unification problem $(P_C)$ has a finite $C$-unification basis; $E$ has \textit{unitary scr-unification type} iff every $C$-unification problem $(P_C)$ has a $C$-mgu.
\end{definition}

When $C$ is empty we have as special case the (standard) notion of \emph{unifier}, \emph{mgu} and \emph{finitary/unitary unification type}; in such case, we indicate unification problems $(P_\emptyset)$ with $(P)$ and the corresponding set of unifiers as $U_E(P)$ instead of  $U_{E}^{svr}(P_\emptyset)$.
%
 The next Proposition (which is an immediate consequence of the above definition and
 Lemma~\ref{Internal Characterization of Admissibility})
 shows the  connection between finitary $C$-unification type and admissibility of $\Pi_2$-rules:

\begin{proposition}\label{Connection between unification and admissibility}
     Assume that $\vdash_{S}$ is decidable. Then if $E$ has finitary svr-unification type (with computable finite $C$-unification bases), then the $\vdash$-admissibility problem for $\Pi_{2}$-rules is decidable too.
 \end{proposition}
%
%


We mention that there is a connection in the reverse direction
that works in case our logic $\vdash$ is decidable, consistent and has a $\bot$-proposition.\footnote{ This means that  there is a constant $\bot\in \fancyL$ such that $(\{\bot\}, \phi)$ belongs to $\vdash$ for all $\phi$.} In fact, in this case $C$-unifiability of $(P_C)$ is equivalent to the non-admissibility of the rule $\forall\overline{p}\Gamma/^{2}\bot$, where $C=\overline{p}$, and where $\Gamma$ is the finite set of formulas obtained by appying the transformers to the equations in $\phi_i=\psi_i$ ($i=1,\dots,k$).

\subsection{Algebraic Characterization of $C$-Unification}\label{subsec:algcharunif}

It will be convenient for our purposes to see unification problems with simple variable restriction from the point of view of finitely presented algebras (following the approach of~\cite{Ghilardi2000} for standard unification problems). For that purpose,
given our equational theory $E$, we will work in $\fpalg{E}$, the opposite of the category of finitely presented $E$-algebras (recall that an algebra is finitely presented iff it is isomorphic to a finitely generated free algebra divided by a finitely generated congruence). Given a finitely presented $E$-algebra $\Alg{A}$, we write it as $\dualg{A}$
when we see it in the opposite category  $\fpalg{E}$ (thus,  $\dualg{A}$ is just a formal
dual of $\Alg{A}$); in particular $\free{X}^{*}$ is the formal dual of $\free{X}$, the free algebra on the finitely many generators $X$. A similar notation is used for morphisms (notice that we have $(\sigma\circ \tau)^{*}= \tau^{*}\circ\sigma^{*}$ for contravariancy). Given an object  $\dualg{B}$ in $\fpalg{E}$, we write $\mathsf{Sub}_{r}(\dualg{B})$ for the set of regular subobjects\footnote{Recall that a regular subobject of $X$ is an equivalence class (wrt to iso) of monomorphisms (with codomain $X$) which happens to be equalizers of a pair of parallel arrows.} of $\dualg{B}$ in the category $\fpalg{E}$; we recall that such regular subobjects correspond, dually, to the finitely presented quotients of $\Alg{B}$ (see e.g. \cite{Ghilardizawadowskybook} for details).

In this context, an $E$-unification problem with simple variable restrictions is a pair $(\Alg{A},C)$, where $C$ is a finite set of free constants, and $\dualg{A}\in \mathsf{Sub}_{r}(\free{X}^{*}\times \free{C}^{*})$. A solution to this problem, which we call suggestively a $C$-\textit{unifier}, is a homomorphism $\sigma:\free{X}\to \free{Z}$, such that $\sigma^*\times 1$ factors in such a way that the diagram of Figure \ref{fig:solutiontounificationproblemwithsimpleconstants} commutes:



\begin{figure}[h]
    \centering
\begin{tikzcd}
\mathbf{F}(X)^{*}\times\mathbf{F}(C)^{*} & \mathbf{F}(Z)^{*}\times\mathbf{F}(C)^{*} \arrow[l, "\sigma^{*}\times 1"'] \arrow[ld, dashed] \\
\mathcal{A}^{*} \arrow[u, " ", tail]         &
\end{tikzcd}    \caption{Solution to Unification Problem with Simple Variable Restrictions}
    \label{fig:solutiontounificationproblemwithsimpleconstants}
\end{figure}

Given such an $E$-unification problem with simple variable restrictions, $(\Alg{A},C)$, and two $C$-unifiers $\sigma:\mathbf{F}(X)\to \mathbf{F}(Z)$ and $\gamma:\mathbf{F}(X)\to \mathbf{F}(W)$, we say that $\sigma$ is more general than $\gamma$, and write $\sigma\leq \gamma$ if there is a homomorphism $k:\mathbf{F}(Z)\to \mathbf{F}(W)$ such that $k\circ \sigma=\gamma$;
as a consequence, the outer triangle of  the following diagram commutes:
\vspace{3mm}

\begin{figure}[h]
\centering
\begin{tikzcd}
& \mathbf{F}(X)^{*}\times\mathbf{F}(C)^{*} &                                                                                       \\
& \mathcal{A}^{*} \arrow[u, " ", tail]         &                                                                                       \\
\mathbf{F}(W)^{*}\times\mathbf{F}(C)^{*} \arrow[ruu, "\gamma^{*}\times1"] \arrow[ru] \arrow[rr, "k^{*}\times 1"'] &                                          & \mathbf{F}(Z)^{*}\times\mathbf{F}(C)^{*} \arrow[luu, "\sigma^{*}\times1"'] \arrow[lu]
\end{tikzcd}
\label{fig:diagramofsection3.1}
\caption{}
\end{figure}
\noindent
(the commutativity of the inner triangle follows as a consequence from the fact that $\mathcal{A}^{*}\hookrightarrow
\mathbf{F}(X)^{*}\times\mathbf{F}(C)^{*}$ is mono).
\begin{remark}\label{rem}
Using the fact that our language $\fancyL$ contains at least a constant symbol,
\emph{we have that  $\sigma\leq \gamma$ iff there is
%
$l:\mathbf{F}(Z)+\mathbf{F}(C)\to \mathbf{F}(W)+\mathbf{F}(C)$ such that $(\sigma^{*}\times 1)\circ l^*=
\gamma^{*}\times 1$}\footnote{
Here $\mathbf{F}(Z)+\mathbf{F}(C)$
denotes the coproduct of
$\mathbf{F}(Z)$ and $\mathbf{F}(C)$
in the category of finitely presented algebras.}. In fact,
homsets among free algebras are not empty, so
if we have  $(\sigma^{*}\times 1)\circ l^*=
\gamma^{*}\times 1$, then letting $l^{*}=\langle l_1^{*},
l_2^{*}\rangle$, we can put $k^*:= l_1^*\circ \langle 1, \alpha^*\rangle$ (where $\alpha$ is any morphism $\mathbf{F}(C)\longrightarrow \mathbf{F}(W)$) and then prove
  $\sigma^{*}\circ k^*= \gamma^{*}$  via elementary properties of products. The latter is seen as follows: from
  $(\sigma^{*}\times 1)\circ l^*=
\gamma^{*}\times 1$, taking first components of pairs,  we get $\sigma^*\circ l_1^*= \gamma^*\circ\pi_{\mathbf{F}(X)}$, so that
$$
\sigma^{*}\circ k^*= \sigma^{*}\circ l_1^*\circ \langle 1, \alpha^*\rangle
=\gamma^*\circ\pi_{\mathbf{F}(X)}\circ \langle 1, \alpha^*\rangle= \gamma^*~~.
$$

\end{remark}

It is easy to see that the above
definition of comparison
for $C$-unifiers gives  a preorder;  we can write $U^{svr}_{E}(\Alg{A},C)$ for the preordered set of $C$-unifiers for $\Alg{A}$.
This `algebraic' approach to $C$-unification is equivalent to the `symbolic' approach of
Definition~\ref{def:utype}, as the following proposition, proved in the Appendix shows:

\begin{proposition}\label{prop:algsymb}
Let  $(P_C)$ a $E$-unification problem with simple variable restriction. If $\Alg{A}$ is a finitely presented algebra with presentation $(P_C)$, then the
antisymmetric quotients of the preordered sets $U^{svr}_{E}(\Alg{A},C)$ and $U_{E}^{svr}(P_C)$
are isomorphic.
\end{proposition}



\section{Interpolation and Finitely Presented Algebras}\label{Algebraic Correspondents to Interpolation}

We analyzed $C$-unification inside the opposite of the category of finitely presented algebras; we now do the same for the interpolation properties mentioned in Theorem \ref{Key Theorem}. We first recall some well-known results
from universal algebra.




\begin{definition}
    We say that a class of $\mathcal{L}$-algebras $\mathbf{K}$ enjoys the property that ``Injections are Transferable" (IT) if whenever $f:\Alg{A}\to \Alg{B}$ is a homomorphism, and $g:\Alg{A}\to \Alg{C}$ is a monomorphism, then there are a morphism $h:\Alg{C}\to \Alg{E}$ and a monomorphism $h':\Alg{B}\to \Alg{E}$ such that $h'\circ f=h\circ g$ (see Figure \ref{fig:itproperty}).
\end{definition}

\begin{figure}[h]
    \centering
\begin{tikzcd}
\mathcal{A} \arrow[r, "f"] \arrow[d, "g"', hook] & \mathcal{B} \arrow[d, "h'", hook] \\
\mathcal{C} \arrow[r, "h"']                      & \mathcal{E}                      
\end{tikzcd}    \caption{Injections are Transferable}
    \label{fig:itproperty}
\end{figure}

The following result is proved in \cite[Lemma 26]{Metcalfe2014}, see also \cite{Czelakowski1985}:

 \begin{theorem}\label{Equivalence of CEP and AP and IT}
     The following are equivalent:
     \begin{enumerate}
         \item $\vdash$ has the Maehara Interpolation Property;
         \item $\mathbf{K}$ has the property (IT).
     \end{enumerate}
 \end{theorem}

We note that it can be shown that (IT) is also equivalent to the conjunction of the congruence extension property and the amalgamation property (this is likewise proven in \cite[Theorem 29]{Metcalfe2014}); we shall make brief use of this further characterization in some examples below. Also notice that  if (IT) holds in $\bf K$, then it holds also in the full subcategory of $\bf K$ formed by the finitely presented algebras: this is because such a subcategory is closed under pushouts, and because a homomorphism $g$ which fits into a factoring of $h=f\circ g$, where $h$ is a monomorphism,
is
itself a monomorphism (so that (IT) for finitely presented algebras comes from the universal property of pushouts).

Right uniform deductive interpolation can likewise be given concrete meaning, in a way which is directly inside finitely presented algebras:

\begin{definition}
    We say that $\mathbf{K}$ is \textit{coherent} if finitely generated subalgebras of finitely presented algebras are again finitely presented.
\end{definition}

The following is shown in \cite[Theorem 2.3]{Kowalski2019}

\begin{theorem}
    $\vdash$ has right uniform deductive interpolation if and only if $\mathbf{K}$ is coherent.
\end{theorem}

The notion of an r-regular category~\cite{Ghilardizawadowskybook} can be useful to summarize the properties of the opposite of the category of finitely presented algebras coming from the above facts.

\begin{definition}
    Let $\mathbf{C}$ be a category. We say that $\mathbf{C}$ is $r$-\textit{regular} if it satisfies the following:
    \begin{enumerate}
        \item it has all finite limits;
        \item epimorphisms are stable under pullback;
        \item every arrow has an epi/regular mono factorization.
    \end{enumerate}
\end{definition}

The next proposition is folklore; it shows the importance of our base assumptions on $E$ (we provide a proof in the appendix, for completeness):

\begin{proposition}
    If $\mathsf{Alg}(E)$ satisfies  (IT) and coherence, then  $\mathsf{Alg}^{op}_{fp}(E)$ is an $r$-regular category.
\end{proposition}

The special property of left-finitary uniform deductive interpolation can likewise be given a straightforward algebraic interpretation inside the category  of finitely presented algebras.
Our notation for r-regular categories is mostly consistent with~\cite{Ghilardizawadowskybook}. In particular, for an object $X$ and an arrow $f:Y\longrightarrow X$, we indicate with $\mathsf{Sub}_{r}(X)$ the poset of the regular subobjects of $X$
and by $f^{-1}: \mathsf{Sub}_{r}(X)\longrightarrow \mathsf{Sub}_{r}(Y)$ the operation of taking pullback along $f$. Projections like $ X\times Y \longrightarrow Y$ are indicated as $\pi_Y$; the identity arrow for an object $X$ (= the maximum regular subobject of $X$) is indicated both with $1_X$ or just with $1$ for simplicity.



\begin{definition}
    Let $\mathbf{C}$ be a category with finite limits. Given $T\in\mathsf{Sub}_{r}(X\times Z)$, we say that a finite collection $B_{1},...,B_{n}\in \mathsf{Sub}_{r}(Z)$ is a $\forall_X$-\textit{factorization} of $T$ (or just a $\forall$-factorization of $T$) if:
    \begin{enumerate}
        \item $\pi_{Z}^{-1}(B_{i})\leq T$ for each $i\leq n$;
        \item for every $C\in\mathsf{Sub}_{r}(Z)$, such that $\pi_{Z}^{-1}(C)\leq T$, there is some $i\leq n$ such that $C\leq B_{i}$.
    \end{enumerate}
    If $n=1$ we say that this is a \textit{singular} $\forall$-factorization. We say that
    $\mathbf{C}$ has the $\forall$-\textit{factorization property} (singular $\forall$-factorization property) if for all objects $X,Z$ and any $S\in \mathsf{Sub}_{r}( X\times Z)$, there is a $\forall_{X}$-factorization (resp. singular $\forall_{X}$-factorization) of $S$.
    The equational theory $E$  (or the equivalent algebraic semantics $\bf K$) has the $\forall$-\textit{factorization property} (singular $\forall$-factorization property) iff
    so does $\mathsf{Alg}_{fp}^{op}(E)$.
\end{definition}

The following  proposition is immediate from the definitions:

\begin{proposition}\label{prop:ffact}
    The logic $\vdash$ has left-finitary deductive interpolation if and only if $\mathsf{Alg}_{fp}^{op}(E)$ has the $\forall$-factorization property.
\end{proposition}

The  $\forall$-\textit{factorization property}  may follow from some natural assumptions: for instance, if $E$ is locally finite it trivially  holds, since there are  only finitely many subobjects. In addition,
 if $E$ is an equational theory axiomatising a logical calculus $\mathsf{L}$, and $\mathsf{L}$ has global post-interpolants, then the $\forall$-factorization property  holds
 and moreover $\forall$-factorizations are singular.

\section{Finitary svr-Unification Types}\label{Unification in Exceptional Logical Systems}

In this section we present a proof of  Theorem \ref{Key Theorem}.
Throughout this section, \emph{we assume that $\mathbf{K}$ 
has (IT), Coherence properties and the $\forall$-factorization property.} 
%
%
%
We  recall the following fact, proved in \cite[Proposition 3.1, pp.51]{Ghilardizawadowskybook}:

\begin{proposition}\label{Existence of Left adjoints with Beck-Chevalley}
    Let $\mathbf{C}$ be an $r$-regular category. Then the pullback functors on regular subobjects have left adjoints satisfying the Beck-Chevalley condition: for every arrow $f:Y\to X$ in $\mathbf{C}$, and every regular subobject $S\in \mathsf{Sub}_{r}(Y)$, there is a regular subobject $\exists_{f}(S)\in \mathsf{Sub}_{r}(X)$ such that:
    \begin{equation*}
        \exists_{f}(S)\leq T \text{ iff } S\leq f^{-1}(T)
    \end{equation*}
    holds for every regular subobject $T\in \mathsf{Sub}_{r}(X)$; in addition, for every pullback square as in Figure \ref{fig:beckchevalleypb}, and every regular subobject $S\in \mathsf{Sub}_{r}(Y_{1})$, the following  condition holds:

    \begin{equation*}
        f_{2}^{-1}(\exists_{f_{1}}(S))=\exists_{p_{2}}(p_{1}^{-1}(S)).
    \end{equation*}

    \begin{figure}[h]
        \centering
\begin{tikzcd}
Z \arrow[d, "p_{2}"'] \arrow[r, "p_{1}"] & Y_{1} \arrow[d, "f_{1}"] \\
Y_{2} \arrow[r, "f_{2}"']                & X
\end{tikzcd}        \caption{Beck-Chevalley Pullback Square}
        \label{fig:beckchevalleypb}
    \end{figure}
\end{proposition}

We now proceed to show that in $r$-regular categories,  $\forall$-factorizations are stable under pullbacks.

\begin{lemma}\label{Stability of forall-factorizations} Let $\mathbf{C}$ be an $r$-regular category;
    consider the pullback of  Figure \ref{fig:pullbackstability}.
    \begin{figure}[h]
        \centering
\begin{tikzcd}
Y\times Z \arrow[r, "f\times 1"] \arrow[d, "\pi_{Y}"'] & X\times Z \arrow[d, "\pi_{X}"] \\
Y \arrow[r, "f"']                                              & X
\end{tikzcd}
\caption{Pullback Stability}
\label{fig:pullbackstability}
    \end{figure}
    If $S\in \mathsf{Sub}_{r}(X\times Z)$ and ${B}_{1},...,{B}_{n}$ is a $\forall_{Z}$-factorization of $S$, then $f^{-1}(B_1),...,f^{-1}({B}_{n})$ is a $\forall_{Z}$-factorization of $(f\times 1)^{-1}(S)$.
\end{lemma}
\begin{proof}
    Let $i\leq n$ be arbitrary, and assume that ${B}_{1},...,{B}_{n}\in \mathsf{Sub}_{r}(X)$ are a $\forall_{Z}$-factorization of $S$. First we want to show that
    \begin{equation*}
        \pi_{Y}^{-1}(f^{-1}({B}_{i}))\leq (f\times 1)^{-1}(S).
    \end{equation*}
    Note that since $\pi_{X}^{-1}({B}_{i})\leq S$, then we have $(f\times 1)^{-1}(\pi_{X}^{-1}({B}_{i}))\leq (f\times 1)^{-1}(S)$ by usual facts on pullbacks. Moreover, since $f\circ\pi_{Y}=\pi_{X}\circ (f\times 1)$ we have that
    $\pi_{Y}^{-1}(f^{-1}({B}_{i}))\leq (f\times 1)^{-1}(S)$.

    For the second property, consider an arbitrary ${C}\in \mathsf{Sub}_{r}(Y)$ such that $\pi_{Y}^{-1}({C})\leq (f\times 1)^{-1}(S)$. Using Proposition \ref{Existence of Left adjoints with Beck-Chevalley}, we can
    apply the functor $\exists_{f\times 1}$, to obtain
    \begin{equation*}
        \exists_{f\times 1}(\pi_{Y}^{-1}({C}))\leq \exists_{f\times 1}(f\times 1)^{-1}(S) \leq S.
    \end{equation*}
    Using Beck-Chevalley, this implies that $\exists_{f\times 1}(\pi_{Y}^{-1}({C}))=\pi_{X}^{-1}(\exists_f({C}))\leq S$. So by the $\forall_{Z}$-factorization property, there is some $B_{i}$ such that $\exists_{f}({C})\leq {B}_{i}$. By  adjunction,  ${C}\leq f^{-1}({B}_{i})$ follows, which shows the property.
\end{proof}

\begin{lemma}\label{Transforming unifiers into inclusions of regular subobjects}
    Given $\dualg{A}\in \mathsf{Sub}_{r}(\mathbf{F}(X)^{*}\times \mathbf{F}(C)^{*})$
    and $\overline{\sigma}:\mathbf{F}(X)\to \mathbf{F}(X')$,
     we have that
    $\overline{\sigma}$
    is a scr-unifier of $(\Alg{A},C)$
     iff
    $1
    \leq (\overline{\sigma}^{*}\times 1)^{-1}(\dualg{A})$.
\end{lemma}
\begin{proof}
    Notice that if $1\leq (\overline{\sigma}^{*}\times 1)^{-1}(\dualg{A})$, then there is an arrow making the triangle of Figure~1 commute. Conversely, if there is such an arrow, then by the universal property of pullbacks if follows that the monic
    $(\overline{\sigma}^{*}\times 1)^{-1}(\dualg{A})\hookrightarrow \mathbf{F}(X')^{*}\times \mathbf{F}(C)^{*}$ has a left inverse, so it is an isomorphism.
\end{proof}

We now prove the key technical result of this section: we  reduce $E$-unification with simple variable restriction to standard $E$-unification:

\begin{theorem}\label{Theorem on Factorization Properties and their connection}
     Let $\Alg{A}^{*}\in \mathsf{Sub}_{r}(\mathbf{F}(X)^{*}\times \mathbf{F}(C)^{*})$, and $\Alg{B}^{*}_{1},...,\Alg{B}^{*}_{n}$ be a $\forall_{\mathbf{F}(C)^{*}}$-factorization of $\Alg{A}^{*}$. We have
    \begin{equation*}
        U^{svr}_{E}(\Alg{A},C)  ~\simeq~ U_{E}(\Alg{B}_{1})\cup\cdots \cup U_{E}(\Alg{B}_{n})
    \end{equation*}
(the isomorphisms being a preordered sets isomorphism).
\end{theorem}
\begin{proof}
    Let $\overline{\sigma}\in U^{svr}_{E}(\Alg{A},C)$ be a scr-unifier; then, by Lemma \ref{Transforming unifiers into inclusions of regular subobjects},
    we have
    $\pi^{-1}_{\mathbf{F}(X')^*}(1)=1\leq (\overline{\sigma}^{*}\times 1)^{-1}(\dualg{A})$.
     Using Lemma \ref{Stability of forall-factorizations} we have that $(\overline{\sigma}^*)^{-1}(\dualg{B}_{1}),...,(\overline{\sigma}^{*})^{-1}(\dualg{B}_{n})$ is a $\forall$-factorization of
     $(\overline{\sigma}^{*}\times 1)^{-1}(\dualg{A})$, hence we have
     $$
     1\leq
     (\overline{\sigma}^{*})^{-1}(\dualg{B}_{i})
     $$
     for some $i$, which means (again by Lemma \ref{Transforming unifiers into inclusions of regular subobjects} applied to the case $C=\emptyset$) that $\overline{\sigma}$ is a unifier of $\Alg{B}_i$.
     Thus we have
     $\overline{\sigma}\in\bigcup_i U_{E}(\Alg{B}_{i})$.

Conversely, if $\overline{\sigma}:\mathbf{F}(X)\to \mathbf{F}(X')$ belongs to $U_{E}(\Alg{B}_{i})$, then by Lemma \ref{Transforming unifiers into inclusions of regular subobjects} (case $C=\emptyset$) we have
 $1_{\mathbf{F}(X')^{*}}
 \leq (\overline{\sigma}^{*})^{-1}(\dualg{B}_{i})$. Since $\pi_{\mathbf{F}(X)^{*}}^{-1}(\dualg{B}_{i})\leq \dualg{A}$ we get
 $\pi_{\mathbf{F}(X')^{*}}^{-1}((\overline{\sigma}^{*})^{-1}(\dualg{B}_{i}))=(\overline{\sigma}^{*}\times 1)^{-1}(\pi_{\mathbf{F}(X)^{*}}^{-1}(\dualg{B}_{i}))\leq (\overline{\sigma}^{*}\times 1)^{-1}(\dualg{A})$. Then by transitivity we obtain
\begin{equation*}
    1\cong \pi_{\mathbf{F}(X')^{*}}^{-1}(1_{\mathbf{F}(X')^{*}})\leq\pi_{\mathbf{F}(X')^{*}}^{-1}((\overline{\sigma}^{*})^{-1}(\dualg{B}_{i}))\leq (\overline{\sigma}^*\times 1)^{-1}(\dualg{A}),
\end{equation*}
 which  means that $\overline{\sigma}\in U_{E}^{svr}(\Alg{A})$ by Lemma \ref{Transforming unifiers into inclusions of regular subobjects}.

 Finally, note that the preorder relation is defined in the same way in
  $U^{svr}_{E}(\Alg{A},C)$ and in   $U_{E}(\Alg{B}_{1})\cup...\cup U_{E}(\Alg{B}_{n})$.
\end{proof}

As a consequence, we can supply the proof of Theorem~\ref{Key Theorem}:

\begin{proof}
    By the results from Section \ref{Algebraic Correspondents to Interpolation}, Maehara Interpolation, Right-Uniform Interpolation and Left-Finitary Interpolation correspond, respectively, to (IT), Coherence and the $\forall$-factorization Property. Consider now the $\Pi_{2}$-rule $\forall\overline{r}\Delta(\overline{q},\overline{r})/^{2}\psi(\overline{q})$ and let $\Theta_{1}(\overline{q}), \dots,\Theta_{n}(\overline{q}) $
    be the finite sets of formulas mentioned in Definition~\ref{def:lint} for $\Delta(\overline{q},\overline{r})$.
    What Theorem~\ref{Theorem on Factorization Properties and their connection} says  (applying the relevant transformers from formulas to equations and back) is that for every $\overline{r}$-invariant substitution $\sigma$, we have  $\vdash \sigma(\Delta)$ iff we have $\vdash \sigma(\Theta_i)$ for some $i$.
    Thus, in view of Lemma~\ref{Internal Characterization of Admissibility}, the rule $\forall\overline{r}\Delta/^{2}\psi$ is admissible iff one of the standard rules $\Theta_i/\psi$ is admissible.
\end{proof}

Theorem~\ref{Theorem on Factorization Properties and their connection} has also some important consequences regarding the decidability of the unification problem with simple variable restrictions:

\begin{corollary}\label{Corollary connecting unification types and admissibility problems}
Suppose that $\mathbf{K}$ satisfies (IT), Coherence and
$\forall$-factorization property.
Then:
\begin{enumerate}
    \item If $E$ has finitary unification type, then it has finitary unification type for the problem with simple variable restrictions.
    \item If $E$ has unitary unification type and $\forall$-factorizations are singular, then it has unitary unification type for the problem with simple constant restrictions.
    \item If $E$-unification is decidable and $\forall$-factorizations are computable, then $E$-unification with simple constant restrictions is decidable as well.
\end{enumerate}
\end{corollary}

\section{Applications}\label{Applications and Negative Results}

In this section we supply some example applications. 

\subsection*{Admissibility of $\Pi_{2}$-rules via Unification}

We can use Corollary~\ref{Corollary connecting unification types and admissibility problems}
and Proposition \ref{Connection between unification and admissibility} in order to directly obtain a decision procedure for admissibility of $\Pi_{2}$-rules via unification.
This goes as follows (let us call $\tau$ the structural transformer from formulas to equations which is granted from the algebraizability hypothesis, see~\cite[Definitions 3.11]{Font2016-dk}): given a $\Pi_{2}$-rule $\forall\overline{p}\Gamma/^{2}\phi$, we first compute a basis of $\overline{p}$-unifiers $\sigma_{0},...,\sigma_{n}$ for the svr-unification problem
given by $\{\tau(\psi)\mid \psi \in \Gamma\}$, and for each of these unifiers -- using the decidability of our logic -- we check whether $\vdash \rho(\sigma_{i}(\tau(\phi)))$ holds or not, where $\rho$ is the inverse transformer of $\tau$. The procedure for computing the basis of $\overline{p}$-unifiers
amounts to the following: using left-finitary deductive uniform interpolation, we compute the finitely many ``approximants" of $\Gamma$ with respect to $\overline{p}$; using decidability of unification, we compute for the transformed equations of each such approximant a finite basis of \emph{standard} unifiers.
As shown by Theorem~\ref{Theorem on Factorization Properties and their connection} above,
these correspond to a basis of $\overline{p}$-unifiers of
$\{\tau(\psi)\mid \psi \in \Gamma\}$. In practical cases, there is no need to apply structural transformers (from formulas to equations and back) because many standard unification algorithms in the literature oriented to propositional logics~\cite{Baader2010}   takes as input directly formulas   (not their transformed  equations).
Thus, below we shall directly speak of `unifiers' and of `$C$-unifiers' of a set of formulas $\Gamma$ (meaning with that the `unifiers' and the `$C$-unifiers' of the transformed set of equations $\{\tau(\psi)\mid \psi \in \Gamma\}$).

As an illustration of how to use our techniques to study admissibility, we turn to the Takeuti-Titani rule, mentioned in Example \ref{Takeuti-Titani Rule}. Such a rule was proven to be admissible over a large class of algebraic signatures, through syntactic methods, by Metcalfe and Montagna \cite{metcalfemontagnasubstructuralfuzzy}, generalising a proof of Baaz and Veith \cite{Baaz2017}.

\begin{example}\label{Example of Takeuti-Titani Rule}
    Let $E$ be the theory of G\"{o}del algebras, i.e., Heyting algebras $\Alg{H}$ satisfying the additional axiom
    \begin{equation*}
        (p\rightarrow q)\vee (q\rightarrow p) ~=~\top.
    \end{equation*}
    Their associated logical system is often denoted by $\mathsf{LC}$ (for `linear calculus'). We will show that the $TT$ rule
    \begin{prooftree}
        \AxiomC{$\forall r\,(g\rightarrow ((p\rightarrow r)\vee (r\rightarrow q)\vee c))$}
        \UnaryInfC{$g\rightarrow (p\rightarrow q) \vee c$}
    \end{prooftree}
is admissible. Following our remarks above, it suffices to show that all of the $C$-unifiers (where $C=\{r\}$) of the the formula $g\to((p\rightarrow r)\vee (r\rightarrow q)\vee c)$ are $C$-unifiers for $g\to((p\rightarrow q)\vee c)$. Since standard $E$-unification is unitary~\cite{wronski} and the conditions of the previous section are satisfied -- we are in the locally finite case and, indeed,  uniform post-interpolants exist -- it suffices to show that the most general
standard unifier for the formula $$\forall_{r}\,(g\to (p\rightarrow r)\vee (r\rightarrow q)\vee c)$$
(namely for the
the  uniform post-interpolant
of $g\rightarrow ((p\rightarrow r)\vee (r\rightarrow q)\vee c)$ wrt $r$)
is also a standard unifier for the consequent $g\rightarrow (p\rightarrow q) \vee c$.

In $\mathsf{LC}$ and in intuitionistic logic systems, uniform interpolants admit a bisimulation semantics which works for Kripke frames, as we proceed to explain. Such semantics can be used to check that a certain formula is really the uniform interpolant of another given one.
The finite Kripke frames corresponding to finite G\"{o}del algebras are precisely the finite frames $\mathfrak{F}=(W,R)$ which are prelinear,
i.e. such that for each $x\in W$, $R[x]=\{y\mid x Ry\}$ is a linear order. By the results from \cite{Ghilardizawadowskybook}, for any prelinear finite Kripke frame, and formula $\phi(\overline{p},\overline{q})$, for each Kripke model $V$ over $\mathfrak{F}$ and over the propositional letters $\overline{p}$, we have for each $x\in W$
\begin{equation*}
    (\mathfrak{F},V),x\sat \forall_{\overline{q}}\phi(\overline{p},\overline{q}) \iff   \text{ for any $(\mathfrak{F}',V')$  $\overline{p}$-bisimilar model}~ (\mathfrak{F}',V'),x\sat \phi(\overline{q},\overline{p}).
\end{equation*}

Using this semantics, we can then show that:
\begin{equation*}
    \forall_{r}\, (g\to (p\rightarrow r)\vee (r\rightarrow q)\vee c) \equiv g\to(p\rightarrow q) \vee c.
\end{equation*}
The right to left side follows from second order intuitionistic propositional logic with the help of the $\mathsf{LC}$-valid formula $(r\to q) \vee (q\to r)$. For the other side we need bisimulation semantics.
Suppose that we have points $x\geq y\in \mathfrak{F}$ such that $(\mathfrak{F},V),x\sat g$,
$(\mathfrak{F},V),x\not\sat c$, $(\mathfrak{F},V),y\sat p$,
$(\mathfrak{F},V),y\not\sat q$.
Form the bisimulation expansion containing a duplicate $y'$ as an immediate successor of $y$  where $y$ refutes $r$, and $y'$ forces $r$. This expansion provides a bisimilar model
such that $g\to (p\rightarrow r)\vee (r\rightarrow q)\vee c$ fails at $x$. Having such an equivalence, the result immediately follows, since
the uniform interpolant we obtained is
precisely the formula in the consequent of the Takeuti-Titani rule.\demo
\end{example}

\subsection*{Finitarity and Unitarity of Unification Type}

In Section \ref{Unification in Exceptional Logical Systems}, we noted that there is a clear connection between svr-unification types and classical unification. It is natural to ask whether in fact the type is always preserved. The next example
shows that it may happen that  $E$-unification is unitary and  svr-unification type is only finitary.

\begin{example}\label{ex:isl}
    Consider the equational theory of \textit{implicative semilattices}, denoted $\mathsf{ISL}$; this corresponds to the $(\top, \wedge,\rightarrow)$-fragment of $\mathsf{IPC}$. Such an equational theory is locally finite, has the amalgamation property and the congruence extension property
    (hence
    it fulfills the hypotheses of Theorem~\ref{Key Theorem}, by the remark we made after Theorem \ref{Equivalence of CEP and AP and IT}).
    It is known that $\mathsf{ISL}$ has unitary elementary unification type \cite{ghilardiunificationprojectivity}. We show that in the setting of unification with simple variable restrictions the unification type becomes finitary. For consider the following svr-unification problem, where we put $C=\{z\}$ (the formula we use is taken from \cite[Example 4.5]{VANGOOLmetcalfetsinakisuniforminterpolation}):
\begin{equation*}
((x\rightarrow z)\wedge (y\rightarrow z)\rightarrow z,\top).
\end{equation*}
    Note that this problem has two incomparable $C$-unifiers, namely $\sigma=\{x\mapsto \top,y\mapsto y\}$ and $\tau=\{y\mapsto \top,x\mapsto x\}$. But we claim that there can be no unifier more general than both of them. For suppose that there was one, say $\mu$. We must have $\mu(x)\not =_E \top$ and $\mu(y)\not =_E \top$  (otherwise $\mu$ would be less general than $\sigma$ or $\tau$). This implies that there are rooted Kriple models (with respective roots $r_1,r_2$) such that
    $(\mathfrak{F}_1,V_1),r_1\not \sat \mu(x)$ and
    $(\mathfrak{F}_2,V_2),r_2\not \sat \mu(y)$. Since $z$ does not occur in $\mu(x), \mu(y)$, we can freely suppose that $(\mathfrak{F}_1,V_1),r_1 \sat z$ and
    $(\mathfrak{F}_2,V_2),r_2 \sat z$. Now build another rooted Kripke model
    $(\mathfrak{F},V)$ by taking the disjoint union of $(\mathfrak{F}_1,V_1)$ and
    $(\mathfrak{F}_2,V_2)$ and by attaching it a new root $r$; we also stipulate that
    $(\mathfrak{F},V),r \not\sat z$. Now then we have that,
    \begin{equation*}
        (\mathfrak{F},V),r \not\sat (\mu(x)\rightarrow z)\wedge (\mu(y)\rightarrow z)\rightarrow z~,
    \end{equation*}
    so $\mu$ cannot be a $C$-unifier\footnote{Incidentally, we notice that the above argument  independently proves that $\mathsf{ISL}$ does not have left uniform interpolation (if it had, by Theorem~\ref{Theorem on Factorization Properties and their connection}, it would also have unitary svr-unification type as it has unitary unification type).}.\demo
\end{example}





\section{svr-Unification in Nuclear Implicative Semilattices}\label{NIS}

As a further
nontrivial
example, we show that the $\{\ell,\wedge,\top, \rightarrow\}$-fragment of lax logic satisfies the hypotheses of Theorem~\ref{Key Theorem}. The  Maehara Interpolation Property follows by the deduction theorem and by inspecting the  proof
of the interpolation property for lax logic in~\cite{rosalie}; Right-Uniform Interpolation and Left-Finitary Uniform Interpolation follow from local finiteness of this variety, shown in ~\cite{BBCGGJ}. The decidability of the admissibility problems for standard rules comes from \emph{finitarity of unification and computability of finite unification bases}:
we will show such properties below (using methods   different methods from those adopted for lax logic in~\cite{GhLe}).

First we need the following folklore fact (implicit in~\cite{G04}):

\begin{proposition}\label{A2}
 Let $\bf K$ be a locally finite variety such that subalgebras of finite projective $\bf K$-algebras are projective. Then unification in $\bf K$ is finitary. Moreover, if $(P)$ is a unification problem, then the unifiers from a finite unification basis for this problem can be chosen so as not to contain more variables than those already occurring in $(P)$.
\end{proposition}

\begin{proof}
 Let $\Alg{A}$ be a finitely presented $\bf K$-algebra, which is finite by local finiteness; let it be a quotient
 $q:\free{X}\longrightarrow \Alg{A}$ of the finitely generated free algebra $\free{X}$.
 Let $\sigma: \free{X} \longrightarrow \free{Y}$ be a unifier for $\Alg{A}$;\footnote{We adopt notation and definitions consistent with those introduced in Subsection~\ref{subsec:algcharunif}. Adopting the approach of~\cite{ghilardiunificationprojectivity,G04} (which views a unifier of $\Alg{A}$ directly as a morphism with domain $\Alg{A}$ and codomain a finitely presented projective algebra) would simplify the arguments.} as such, $\sigma$ factors through
$q$ as $\overline{\sigma}\circ q$. Taking the image factorization, we can further factorize $\overline{\sigma}$ as $\Alg{A}\buildrel{\overline{\sigma}_0} \over \longrightarrow P\buildrel{\iota} \over\longrightarrow \free{Y}$, where $P$ is projective as a subalgebra of a projective algebra. By projectivity, the surjective map $\overline{\sigma}_0\circ q:\free{X}\longrightarrow P$ has a section $s$ (i.e. $\overline{\sigma}_0\circ q\circ s=1_{P}$). Now
$s\circ \overline{\sigma}_0\circ q: \free{X}\longrightarrow \free{X}$ is another unifier for $\Alg{A}$ (because it factors though $q$) and is more general than $\sigma$ because
$$
\sigma\circ (s\circ \overline{\sigma}_0\circ q)= (\iota\circ \overline{\sigma}_0\circ q) \circ (s\circ \overline{\sigma}_0\circ q)= \iota\circ \overline{\sigma}_0\circ q=\sigma
$$
Thus unifiers of $\Alg{A}$ having domain and codomain $\free{X}$ form a unification basis. Since $\mathbf{F}(X)$ is finite, and hence there can be only finitely many unifiers with this domain and codomain, this unification basis is finite.
\end{proof}

The variety algebraizing the $\{\ell,\wedge,\top, \rightarrow\}$-fragment of lax logic is formed by \emph{nuclear implicative semilattices}, namely by the algebras $\Alg{A}= (A, \top, \wedge, \rightarrow, \ell)$,  where $(A, \top, \wedge, \rightarrow)$ is an implicative semilattice and $\ell: A\longrightarrow A$ is a \textit{nucleus}, i.e., a unary operator satisfying the conditions
$$
x\leq \ell x, ~~~~~\ell(x\wedge y)=\ell x \wedge \ell y, ~~~~~ \ell \ell x \leq \ell x~~.
$$
The category of finite nuclear semilattices ${\bf NIS}_{fin}$ (which by local finiteness, coincides with the category of finitely presented such algebras) is dual to the category ${\bf SF}_{fin}$ of finite $S$-posets and morphisms~\cite{BBCGGJ}. An $S$-poset is a triple $(X, \leq, S)$, where $(X,\leq )$ is a poset and $S\subseteq X$ is a subset; a morphism
$$f: (X, \leq, S)\longrightarrow (Y, \leq, T)$$ between $S$-posets is a partial map $f$ satisfying the following conditions
(we let
$dom(f)$ be the domain of $f$, and
$x<y$ mean $x\leq y$ and $x\neq y$):
\begin{description}
 \item[{\rm (i)}] if $x<y$ and $x,y\in dom(f)$ then $f(x)< f(y)$;
 \item[{\rm (ii)}] if $x\in dom(f)$ and $f(x)<y$ there there is $x'$ such that $x<x', x'\in dom(f)$ and $f(x')=y$;
 \item[{\rm (iii)}] $f^{-1}(T)= dom(f)\cap S$;
 \item[{\rm (iv)}] if $s\in S, s\leq x$ and $x\in dom(f)$, then there are $s', x'\in dom(f)$ such that $s\leq s'\leq x'$, $s'\in S$ and $f(x')=f(x)$.
\end{description}
We also need the following fact~\cite[Proposition 5.1]{BBCGGJ}:

\begin{lemma}
 The dual of $f: (X, \leq, S)\longrightarrow (Y, \leq, T)$ is injective iff $f$ is surjective and the dual of $f$ is surjective iff $f$ is injective and totally defined.
\end{lemma}

\vspace{1mm}

Let $\alpha$ be an antichain in an $S$-poset $(X, \leq, S)$, i.e., a set of mutually $\leq$-incomparable elements; a \emph{cover} of $\alpha$ is some $s\in S$ such that
$\alpha$ is the set of the immediate successors of $s$.

\begin{lemma}
 The dual of a finite $S$-poset $(X, \leq, S)$ is a projective nuclear implicative semilattice iff every antichain $\alpha\subseteq X$ such that $\alpha\not \subseteq S$ has a cover.
\end{lemma}

\begin{proof}
 If the dual of $(X, \leq, S)$ is  projective and $\alpha\subseteq X$ is such that $\alpha\not \subseteq S$ in order to find a cover of $\alpha$ it is sufficient to embed $(X, \leq, S)$ into the $S$-poset $(X\cup\{*\}, \leq, S\cup\{*\})$ obtained from $(X, \leq, S)$ by adding an extra element $*$ covering $\alpha$: we show that the retract $r:(X\cup\{*\}, \leq, S\cup\{*\})\longrightarrow (X, \leq, S)$  must map $*$ to a cover $r(*)$ of $\alpha$. First, since $r$ is a retract, we must have $r(x)=x$ for $x\in X$. Then, $*$ must be in the domain of $r$ by (i)-(iv): in fact, we have $*\leq d$ for $d\in\alpha\setminus S$ and so  there must be $s'$ such that $s'\in dom(r)\cap (S\cup \{*\})$ and $*\leq s'\leq d$. This $s'$ must coincide with $*$ (because $d$ is an immediate successor of $*$, and $d\not\in S$), so $*\in dom(r)$; also $r(*)\in S$ by (iii) and $r(*)$ covers $\alpha$ by (i)-(ii).

 Conversely, suppose that $(X, \leq, S)$ satisfies the condition of the lemma and take a total embedding into a finite $S$-poset $ (X, \leq, S)\hookrightarrow  (Y, \leq, T)$. We find a retract $r: (Y, \leq, T) \longrightarrow (X, \leq, S)$ by defining $r(y)$ by induction on the height of $y\in Y$. If $y\in X$, we put $r(y)=y$; otherwise let $\alpha$ be the antichain  of the minimal elements of
 $\{ r(y') \mid y < y' ~\&~y'\in dom(r)\}$.
 If $y\not \in T$ or $\alpha\subseteq S$, then $dom(r)$ will not include $y$, otherwise we take $r(y)$ to be a cover of $\alpha$. It is easy to check that $r$ is a retract of the inclusion $ (X, \leq, S)\hookrightarrow  (Y, \leq, T)$. To prove that $r$ satisfies the above conditions (i)-(iv), one shows by induction on the height of $y$ that such conditions are satisfied by the restriction of $r$ to the cone $\{ y'\in Y\mid y\leq y'\}$.
\end{proof}

\begin{lemma}\label{lem:nisl}
 Subalgebras of projective finite nuclear implicative semilattices are projective.
\end{lemma}

\begin{proof}
 Let $(X, \leq, S)$ be dual to a finite projective nuclear implicative semilattice and let $f: (X, \leq, S)\longrightarrow (Y, \leq, T)$ be surjective. Take an antichain $\alpha\subseteq Y$ such that $\alpha\not\subseteq T$; let $\beta$ be the antichain formed by the minimal elements of $f^{-1}(\alpha)$; we have that $f(\beta)=\alpha$ and $\beta\not\subseteq S$ by condition (iii) above. Thus there is a cover $s\in S$ for $\beta$ and by (iv) we must have that $s\in dom(f)$ and that $f(s)$ is a cover of $\alpha$.
\end{proof}

\begin{theorem}
 Admissibility of standard rules (and consequently also of $\Pi_2$-rules) is decidable in the $\{\ell,\wedge,\top, \rightarrow\}$-fragment of lax logic.
\end{theorem}

\begin{proof}
 By Proposition~\ref{A2} and Lemma~\ref{lem:nisl},
 standard unification for this logic has finitary unification type and finite bases of unifiers are computable.
 This guarantees decidability of admissibility of standard rules and also of $\Pi_2$-rules, as explained at the beginning of Section~\ref{Applications and Negative Results}.
\end{proof}

\section{Conclusions and Further Work}\label{Conclusions and Further Work}

In this work we  analyzed a new type of unification problems which are properly situated in the literature between elementary unification and the so-called unification with linear constant restrictions. Their interest here lies primarily in the connection with admissibility of non-standard $\Pi_{2}$-rules over logics.
We supplied some first results and we discuss here
 several natural  open questions.

One natural question that applies to logics, is whether the admissibility of more complex logical rules could be related to unification with linear constant restrictions.
Such a question can be motivated also from a model theoretic point of view, since it
is connected to decision problems for the positive theory of  free algebraic structures.


Finally, and most importantly, it would be interesting to explore unification with simple variable restrictions for systems  not covered by the results of the present paper. A natural example in this sense is the modal logic system $\mathsf{S4}$, which is well-known not to enjoy uniform interpolation \cite{Ghilardi1995}.
We note that the obvious approach to attack this problem -- a generalization of the \textit{projective approximations} from \cite{Ghilardi2000} -- does not work in the obvious way, since the key technique of Lowenheim Substitutions seems  not to be available.


\Appendix

In this appendix we collect some missing
(mostly folklore)
technical proofs.

\subsection*{(A) Proof of Proposition~\ref{prop:algsymb}}


We recall that, given  finite sets of variables $Z,Z'$ and a substitution $\sigma$ with domain $Z$ and codomain $Z'$, we can canonically associate with it the homomorphism $\eta(\sigma):\mathbf{F}(Z) \longrightarrow \mathbf{F}(Z')$ mapping the equivalence class of a term $t\in Fm_{\fancyL}(Z)$ to the equivalence class of the term $\sigma(t) \in Fm_{\fancyL}(Z')$. This correspondence is \emph{bijective}, if we identify substitutions up to $=_E$: for that reason we always used the same letters $\sigma, \tau, \dots$ for substitutions and free algebra homomorphisms. In this subsection, however, we conveniently distinguish between $\sigma$ and the associated homomorphism $\eta(\sigma)$. Also note that the correspondence $\sigma \mapsto \eta(\sigma)$ commutes with compositions, in the sense that it maps the composition of substitutions into the compositions of homomorphisms.

In addition, notice that if $X,X',C$ are finite disjoint sets, and $\eta(\sigma):\mathbf{F}(X\cup C)\longrightarrow \mathbf{F}(X'\cup C)$ is the free algebras homomorphism induced by  a $C$-invariant substitution $\sigma$,
then there is a substitution $\overline{\sigma}:\mathbf{F}(X)\to \mathbf{F}(X')$
    such that $\eta(\sigma)\simeq\eta(\overline{\sigma})+1$, up to the isomorphisms
    $\mathbf{F}(X)+\mathbf{F}(C)\simeq\mathbf{F}(X\cup C)$ and
    $\mathbf{F}(X')+\mathbf{F}(C) \simeq\mathbf{F}(X'\cup C)$.

\vskip 2mm \noindent
\textbf{Proposition~\ref{prop:algsymb}} \textit{
Let  $(P_C)$ a $E$-unification problem with simple variable restriction. If $\Alg{A}$ is a finitely presented algebra with presentation $(P_C)$, then the
antisymmetric quotients of the preordered sets $U^{svr}_{E}(\Alg{A},C)$ and $U_{E}^{svr}(P_C)$
are isomorphic.}

\begin{proof}
Let $(P_C)$ be the unification problem
    with simple variable restriction
    $$
    (P_C)~=~(s_{1},t_{1}),...,(s_{k},t_{k})
    $$
    where $X=\{x_{1},...,x_{n}\}\cup C=\{c_{1},...,c_{m}\}$ are the disjoint sets of variables occuring in these terms. The corresponding finitely presented algebra algebra $\Alg{A}$ is $\mathbf{F}(X\cup C)/S$ where $S$ is the smallest congruence generated by the set of pairs $\{(s_{1},t_{1}),...,(s_{k},t_{k})\}$.

Define now a map $e:U_{E}^{svr}(P_C)\to U_{E}^{svr}(\Alg{A},C)$ as follows: let $\sigma\in U_{E}^{svr}(P_C)$, and suppose that $Y\cup C$ is its codomain, where $Y$ are the variables
     occurring in the terms $\sigma(x_{1}),...,\sigma(x_{n})$ (note that this set is disjoint from $C$). Thus $\sigma$ can be restricted to a substitution $\overline{\sigma}$ with domain $X$ and codomain $Y$. Such $\sigma$ and $\overline{\sigma}$ induce  homomorphisms
     $$
     \eta(\overline{\sigma}): \mathbf{F}(X)\longrightarrow \mathbf{F}(Y), ~~~~
     \eta({\sigma}): \mathbf{F}(X\cup C)\longrightarrow \mathbf{F}(Y\cup C)
     $$
such that $\eta(\overline{\sigma})+1=\eta({\sigma})$. We put $e(\sigma):=\eta(\overline{\sigma})$.
      Recalling the definitions from Section~\ref{Unification with Simple Constant Restrictions} and Subsection~\ref{subsec:algcharunif}, it is clear that
     $\sigma\in U_{E}^{svr}(P_C)$ iff $e(\sigma)\in U^{svr}_{E}(\Alg{A},C)$ (in fact $e(\sigma)^*\times 1=\eta(\overline{\sigma})^* \times 1=\eta(\sigma)^*$ factors through $\Alg{A}$ iff the kernel of $\eta(\sigma)$ contains the equivalence classes
     in $\mathbf{F}(X\cup C)$
     of the pairs of terms
     $(s_{1},t_{1}),...,(s_{k},t_{k})$, which precisely means that $\sigma$ unifies them).

     The map $e:U_{E}^{svr}(P_C)\to U_{E}^{svr}(\Alg{A})$ is bijective up to $=_E$ equivalence of substitutions, thus it becomes a real bijection when we identify substitutions up to the comparison order. This order is preserved and reflected by $e$ if we compare
    the preordered sets $U^{svr}_{E}(\Alg{A},C)$ and $U_{E}^{svr}(P_C)$  using the equivalent definition for $U^{svr}_{E}(\Alg{A},C)$ given by the Remark~\ref{rem} of Subsection~\ref{subsec:algcharunif}.
    In fact, for two substitutions $\sigma_1: Fm_{\fancyL}(Z)\longrightarrow Fm_{\fancyL}(Z_1)$ and
    $\sigma_2: Fm_{\fancyL}(Z)\longrightarrow Fm_{\fancyL}(Z_2)$, we have that $\sigma_1\leq \sigma_2$ iff there is  a substitution $\theta: Fm_{\fancyL}(Z_2)\longrightarrow Fm_{\fancyL}(Z_1)$ such that $\eta(\theta\circ\sigma_2)=\eta(\theta)\circ \eta(\sigma_2)= \eta(\sigma_1)$.
\end{proof}

\vspace{1mm}

\subsection*{(B) r-Regularity of $\mathsf{Alg}^{op}_{fp}(E)$}

\noindent
\textbf{Proposition A.1} \textit{
    Let $E$ be an equational theory enjoying coherence and (IT). Then $\mathsf{Alg}^{op}_{fp}(E)$ is an $r$-regular category.}
\begin{proof}
    The fact that $\mathsf{Alg}^{op}_{fp}(E)$ has all finite limits is immediate, given the standard fact that
     $\mathsf{Alg}^{op}_{fp}(E)$
     has all finite colimits. To see the factorization properties, assume that $f^{*}:\Alg{B}^{*}\to \Alg{A}^{*}$ is a map of $\mathsf{Alg}^{op}_{fp}(E)$; then $f:\Alg{A}\to \Alg{B}$ is a homomorphism, which as usual has an image factorization
    \begin{equation*}
        \Alg{A}\xrightarrow{f} Im(f) \xrightarrow{i} \Alg{B}.
    \end{equation*}
    Note that $Im(f)$ is a finitely generated subalgebra of $\Alg{B}$, since it is a quotient of the finitely presented algebra $\Alg{A}$; hence by coherence, $Im(f)$ is itself finitely presented; this means that the image factorization lives inside of the category $\mathsf{Alg}^{op}_{fp}(E)$, and hence by duality, $\mathsf{Alg}^{op}_{fp}(E)$ has the desired factorization properties.

    Finally, assume that we have a pullback square as in the Figure below, and that $g^*$ is an epimorphism. By (IT) and duality, there exist an epimorphism $p:\Alg{E}^{*}\to \Alg{B}^{*}$ and a homomorphism $p':\Alg{E}^{*}\to \Alg{C}^{*}$ commuting the outer square below. Since $\Alg{D}^{*}$ is a pullback, there is a connecting morphism $k:\Alg{E}^{*}\to \Alg{D}^{*}$; but this means that $h_{2}$ is an epimorphism as well, since
    it is the second component of
    an epimorphism.
    \end{proof}

    \begin{figure}[h]
        \centering
\begin{tikzcd}
\mathcal{E}^{*} \arrow[rrd, "p", bend left] \arrow[rdd, "p'"', bend right] \arrow[rd, "k", dashed] &                                                        &                           \\
  & \mathcal{D}^{*} \arrow[d, "h_{1}"'] \arrow[r, "h_{2}"] & \mathcal{B}^{*} \arrow[d] \\
& \mathcal{C}^{*} \arrow[r, "g^*"']                      & \mathcal{A}^{*}
\end{tikzcd}

    \end{figure}

\bibliographystyle{aiml22}
\bibliography{aiml22}

\end{document}